\documentclass{amsart}

\usepackage[T1]{fontenc}
\usepackage[utf8]{inputenc}
\usepackage{hyperref}
\usepackage{amsmath}
\usepackage{amssymb}
\usepackage{amscd}
\usepackage{graphicx}
\usepackage{mathdots}
\usepackage{gensymb}
\usepackage{epstopdf}
 \usepackage{todonotes}
 \usepackage{youngtab}
 \usepackage{wrapfig}

\usepackage{epsfig}       %For postscript
\usepackage{epic,eepic}        %For epic and eepic output from xfig

\definecolor{red}{RGB}{191,0,0}

\newtheorem{thm}{Theorem}[section]

\newtheorem{conj}[thm]{Conjecture}
\newtheorem{quest}[thm]{Question}

\theoremstyle{definition}

\theoremstyle{remark}

\numberwithin{equation}{section}

\newcommand{\group}[1]{\mathrm{#1}}

\newcommand{\Cat}{\operatorname{Cat}}

\begin{document}

\title{Procesi's Conjecture on the Formanek-Weingarten Function is False}
%\keywords{Weingarten function, Catalan numbers, join-cut equation.}
%\subjclass[2010]{test}

\author[M.~Dołęga]{Maciej Dołęga}
\address{
Institute of Mathematics, 
Polish Academy of Sciences, 
ul. Śniadeckich 8, 
00-956 Warszawa, Poland.
}
\email{mdolega@impan.pl}

\author{Jonathan Novak}
\address{Department of Mathematics, University of California, San Diego, USA}
\email{jinovak@ucsd.edu}

\thanks{MD is supported from {\it Narodowe Centrum Nauki},
  grant 2021/42/E/ST1/00162/}
 \thanks{JN is supported by NSF 
  grant DMS 1812288 and a Lattimer Fellowship.}

\maketitle

\begin{abstract}
	In this paper, we disprove a recent monotonicity conjecture of C. Procesi on the 
	generating function for monotone walks on the symmetric group, an object which
	is equivalent to the Weingarten function of the unitary group.
\end{abstract}

\section{Introduction}
Let $\Gamma_d$ be the Cayley graph of the symmetric group $\group{S}(d)$ as generated by 
the conjugacy class of transpositions. Thus $\Gamma_d$ is a
${d\choose 2}$-regular graded graph with levels $L_0,\dots,L_{d-1},$
where $L_k$ is the set of permutations which factor into a product of 
$d-k$ disjoint cycles. Let us mark each edge of $\Gamma_d$ corresponding 
to the transposition $(i\ j)$ with $j \in \{2,\dots,d\}$, the larger of the two symbols interchanged. This edge labeling was first considered by Stanley 
\cite{Stanley} and Biane \cite{Biane} in 
connection with noncrossing partitions and parking functions.

A walk on $\Gamma_d$ is said to be \emph{monotone} if the labels of the edges it traverses
form a weakly increasing sequence. 
The combinatorics of such walks has been intensively studied in recent years,
beginning with the discovery \cite{Novak:Banach} that these
trajectories play the 
role of Feynman diagrams for integration with respect to Haar measure on 
unitary groups. This is part of a broader subject nowadays 
known as Weingarten calculus, see \cite{CMN}. 

Although non-obvious, it is a fact that the number of monotone
walks of given length $r$ between two given permutations $\rho,\sigma \in \group{S}(d)$
depends only on the cycle type $\alpha \vdash d$ of the permutation $\rho^{-1}\sigma$. 
It is therefore sufficient to consider the number $m^r(\alpha)$ of $r$-step monotone walks
on $\Gamma_d$ beginning at the identity permutation and ending at a fixed permutation of cycle
type $\alpha$. To each partition $\alpha \vdash d$ we associate the generating function 

	\begin{equation}
		M_\alpha(x) = \sum_{r=0}^\infty m^r(\alpha) x^r
	\end{equation}
	
\noindent
enumerating monotone walks on $\Gamma_d$ of arbitrary length and type $\alpha$. 
It is known \cite{MN} that 

	\begin{equation}
	\label{eqn:CharacterFormula}
		M_\alpha(x) = \sum_{\lambda \vdash d} \frac{\chi^\lambda_\alpha}{\prod_{\Box \in \lambda} h(\Box) (1-c(\Box)x)},
	\end{equation}

\noindent
where $\chi^\lambda_\alpha$ are the irreducible characters of the symmetric 
group $\group{S}(d)$, with $h(\Box)$ and $c(\Box)$ being, respectively, the hook length 
and content of a given cell $\Box$ in the Young diagram of $\lambda$
(see \cite{Stanley:EC2} for definitions). 
In particular, $M_\alpha(x)$ is a rational function of $x$
which may be considered as a continuous function of $x$ on  
the interval $(0,\frac{1}{d-1})$ whose outputs are positive
rational numbers. Up to a simple rescaling, the values 
$M_\alpha(\frac{1}{N})$ coincide with the values
of the Weingarten function of the 
unitary group $\group{U}(N)$; see \cite{MN,Novak:Banach}.

In a recent paper~\cite{Procesi}, Procesi has pointed out that
the function $M_\alpha(x)$ was also studied from the 
perspective of classical
invariant theory by Formanek, and that in this context the 
values $M_\alpha(\frac{1}{d})$ have special significance. 
Procesi tabulated these numbers for all diagrams $\alpha \vdash d \leq 8$,
and on the basis of these computations made the following conjecture.

	\begin{conj}
	\label{conj:Procesi}
		If $\alpha > \beta$ in lexicographic order, then $M_\alpha(\frac{1}{d}) > M_\beta(\frac{1}{d}).$
	\end{conj}
	
In this brief note we give explicit numerical examples which show 
that Conjecture \ref{conj:Procesi} is false.
	
\section{Small $x$}
We first clarify that Procesi's Conjecture \ref{conj:Procesi} refers 
to lexicographic order on
partitions viewed as nondecreasing sequences of positive integers,
with $1$ the first letter in the alphabet, $2$ the second letter,
and so on. For example,
the partitions of six listed in lexicographic order are

	\begin{align*}
		&111111 \\
		&11112 \\
		&1113 \\
		&1122 \\
		&114 \\
		&123 \\
		&15 \\
		&222 \\
		&24 \\
		&33 \\
		&6,
	\end{align*} 

\noindent
and Conjecture \ref{conj:Procesi} says that the numbers $M_\alpha(\frac{1}{6})$ 
strictly decrease as $\alpha$ moves down this list, and this is so.
However, the pattern fails for sufficiently large degree $d$.

The first sign that Conjecture \ref{conj:Procesi} 
might be false in general is that it is incompatible
with the known $x \to 0$ asymptotics of $M_\alpha(x)$. The minimal length of a
walk on $\Gamma_d$ from the identity to a permutation of type $\alpha$ is 
$d-\ell(\alpha)$, and thus by parity the number $m^r(\alpha)$ can only be positive
when $r=d-\ell(\alpha)+2k$ with $k$ a nonnegative integer. 
We may therefore reparameterize  the counts $m^r(\alpha)$
as $m_k(\alpha) := m^{d-\ell(\alpha)+2k}(\alpha)$
for $k \in \mathbb{N}_0$. 
The generating function $M_\alpha(x)$ then becomes

	\begin{equation}
		M_\alpha(x) = x^{d-\ell(\alpha)} \sum_{k=0}^\infty m_k(\alpha) x^{2k}.
	\end{equation}
	
\noindent
It is then clear that 

	\begin{equation}
		\lim_{x \to 0} \frac{M_\beta(x)}{M_\alpha(x)} = 0
	\end{equation}
	
\noindent
whenever $\ell(\alpha)>\ell(\beta)$, which is incompatible with lexicographic order. 

One might nevertheless hope that when we compare the small $x$ behavior of $M_\alpha(x)$ and $M_\beta(x)$ with $\alpha$ and $\beta$ being partitions of the same length, we find compatibility with lexicographic order. This too is false, as can be 
seen from the fact \cite{MN} that 

	\begin{equation}
	\label{eq:Catalan}
		m_0(\alpha) = \prod_{i=1}^{\ell(\alpha)} \Cat_{\alpha_i-1},
	\end{equation}
	
\noindent
where $\Cat_n=\frac{1}{n+1}{2n \choose n}$ is the Catalan number.
Then for $\alpha,\beta \vdash d$ partitions of the same length $\ell$,
we have

	\begin{equation}
		\lim_{x\to 0} \frac{M_\beta(x)}{M_\alpha(x)} =\prod_{i=1}^\ell \frac{\Cat_{\beta_i-1}}{\Cat_{\alpha_i-1}}.
	\end{equation}
	
\noindent
For small values of $d$, it does indeed appear to be the case that
this product is smaller than $1$ when $\alpha > \beta$, but this is 
a law of small numbers. Consider the case 
where

    \begin{equation}
        \alpha = (1,\underbrace{3,\dots,3}_n)
        \quad\text{ and }\quad
        \beta = (\underbrace{2,\dots,2}_{n},n+1)
    \end{equation}
    
\noindent  
Then $\alpha$ and $\beta$ are partitions of the same
degree $d = 3n+1$, they have the same length $\ell(\alpha)=\ell(\beta)=n+1$,
and $\alpha$ precedes $\beta$ in the lexicographic order. 
However, the ratio of the corresponding Catalan products 
tends to infinity as $n \to \infty$, 

    \begin{equation}
        \frac{\Cat_{n}}{2^n} \sim \frac{1}{\sqrt{\pi} n^{3/2}}\cdot 2^n.
    \end{equation}

	\section{Counterexamples}
	To give a counterexample to Conjecture \ref{conj:Procesi} itself,
	we return to the character formula \eqref{eqn:CharacterFormula},
	which in fact yields counterexamples if one goes a bit 
	farther than the data tabulated in \cite{Procesi}.
	Let $\alpha^+$ denote the successor of $\alpha$ in the lexicographic order. The first value of $d$ for which Conjecture~\ref{conj:Procesi} fails is 
	the famously unlucky number $d=13$,
	for which there exists precisely one violating pair $\alpha,\alpha^+$.
	This pair is
	
	\[ M_{(1^6,7)}\left(\frac{1}{13}\right) = \frac{13^{13}}{(13!)^2}\frac{30132115571}{1149266300} <\frac{13^{13}}{(13!)^2}\frac{426729597219}{16089728200} = 
	M_{(1^5,2^4)}\left(\frac{1}{13}\right)\]

        We have tested Conjecture~\ref{conj:Procesi} for $d \leq 20$ and
        it fails for all $13 \leq d \leq 20$. Moreover the size of the
        set
        \[ G_d := \left\{ \alpha \vdash d\colon M_{\alpha}\left(\frac{1}{d}\right)<M_{\alpha^+}\left(\frac{1}{d}\right)\right\}\]
of consecutive failures at rank $d$ 
increases with $d$. For instance
        \begin{align*}
G_{14}= \{&(1^7,7), (1^5,2,7), (1^5,9)\},\\
G_{15} =\{&(1^8,7), (1^6,2,7), (1^6,9), (1^4,11), (1^3,2,10), (1^3,3,9)\},
                  \\
G_{16}=\{&(1^{11},5), (1^9,7), (1^7,2,7), (1^7,9), (1^6,10), (1^5,2^2,7), (1^5,11), (1^4,2,10),\\ &(1^4,3,9), (1^3,13),  (1,4,11)\}.         
        \end{align*}

        Even though Conjecture~\ref{conj:Procesi} seems to fail for all
        $d \geq 13$ the structure of
        the failure set $G_d$ appears to be very interesting: it seems
        that when $d$ is large, 
        the points in $G_d$ form many short 
        lexicographic intervals and one large lexicographic interval. For
        instance $|G_{20}| = 45$, so the proportion of the length of
        a typical interval on which $M_{\alpha}(\frac{1}{|\alpha|})$
        is monotone is equal to $\frac{1}{45}$. Nevertheless, for the
        interval $((1,2^2,4,11), (2,5,13)]$, whose
        cardinality is equal to $151$, one has $((1,2^2,4,11),
        (2,5,13)] \cap G_{20} = \{(2,5,13)\}$. The number of partitions of size $20$ is
        $627$, therefore there exists
        an interval on which $M_{\alpha}(\frac{1}{|\alpha|})$ is
        monotone and which is more than ten times longer than its
        expected length. This suggests that a weaker version of
        Conjecture~\ref{conj:Procesi} might be true.
        Let $\mathcal{P}_d$ denote the set of partitions of size $d$.

        	\begin{quest}
	\label{conj:ProcesiBis}
		Is it true that there exists constant $C>0$ such that
                for every positive integer $d$ there exists partitions
                $\alpha^d > \beta^d \in \mathcal{P}_d$ such that for
                every lexicographic sequence 
                $\alpha^d \geq \alpha > \beta \geq \beta^d $ we have 
                $M_\alpha(\frac{1}{d}) > M_\beta(\frac{1}{d})$ and
                $\frac{|[\alpha_d,\beta_d]|}{|\mathcal{P}_d|} \geq C$?
              \end{quest}

              We do not know the answer to this question and we leave it wide open. It 
              would also be very interesting to find an explicit
              description of the set $G_d$, which appears to consists of 
              very specific partitions which might be classifiable. Even though
              Conjecture~\ref{conj:Procesi} turned out to be false we
              believe that the research initiated by
              Procesi~\cite{Procesi} on the behaviour of the function
              $M_{\alpha}(\frac{1}{|\alpha|})$ merits further
              investigation. Indeed, Procesi's work has added a
              new and largely unexplored dimension to 
              Weingarten calculus.

              \section*{Acknowledgments}
The SageMath computer algebra system
\cite{sagemath} has been used for experimentation leading up to the results
presented in the paper.

\bibliographystyle{amsplain}

\begin{thebibliography}{10}

    \bibitem{Biane}
	P. Biane,
	\textit{Parking functions of types A and B},
	Electron. J. Combin. \textbf{9} (2002), \#N7.
	
	\bibitem{CMN}
	B. Collins, S. Matsumoto, J. Novak,
	\textit{The Weingarten calculus},
	Not. Amer. Math. Soc., in press.

	\bibitem{MN}
	S. Matsumoto, J. Novak,
	\textit{Jucys-Murphy elements and unitary matrix integrals},
	Int. Math. Res. Not. IMRN \textbf{2} (2013), 362-397.
	
	\bibitem{Novak:Banach}
	J. Novak,
	\textit{Jucys-Murphy elements and the Weingarten function},
	Banach Cent. Publ. \textbf{89} (2010), 231-235.
	
	\bibitem{Procesi}
	C. Procesi,
	\textit{A note on the Formanek Weingarten function},
	Note Mat. \textbf{41} (2021), 69-109.

        \bibitem{sagemath}
          The Sage Developers,
          \textit{{S}agemath, the {S}age {M}athematics {S}oftware
  {S}ystem}, {\tt https://www.sagemath.org}.
	
	\bibitem{Stanley:EC2}
	R. P. Stanley,
	\textit{Enumerative Combinatorics}. Vol. 2. Cambridge University Press,
	New York, 1999.
	
	\bibitem{Stanley}
	R. P. Stanley, 
	\textit{Parking functions and noncrossing partitions},
	Electron. J. Combin. \textbf{4} (1997), \#2.

	
\end{thebibliography}

\end{document}